\documentclass[a4paper, 11pt]{article}
\usepackage[english]{babel}
\usepackage{latexsym}
\usepackage{amsfonts, amsthm}
\usepackage{amssymb}
\usepackage{amscd}
\usepackage{amsmath}
\usepackage{graphics}
\usepackage{eepic}
\usepackage{epsfig}
\usepackage[noxcolor]{pstricks}
\newtheorem{theorem}{Theorem}%[section]
\newtheorem{lemma}[theorem]{Lemma}
\newtheorem{proposition}[theorem]{Proposition}

\newtheorem{definition}%[theorem]
{Definition}
\newtheorem{remark}{Remark\/}
\newtheorem{example}{Example\/}
{Conjecture\/}

\newcommand{\C}{{\mathbb C}}
\newcommand{\Z}{{\mathbb Z}}

\newcommand{\I}{\mathcal{I}}

\begin{document}

\title{\textbf{Poincar\'e series of embedded filtrations}}
\author{Ann Lemahieu \footnote{Ann Lemahieu; K.U.Leuven, Departement Wiskunde,
Celestij\-nenlaan 200B, B-3001 Leuven, Belgium, email:
ann.lemahieu@wis.kuleuven.be. The
research was partially supported by the Fund of Scientific Research
- Flanders and MEC PN I+D+I MTM2007-64704.}
\date{}}
\maketitle {\footnotesize \emph{\noindent \textbf{Abstract.---}
In this article we define a Poincar\'e series on a subspace of a complex analytic germ, induced
by a multi-index filtration on the ambient space.
We compute this Poincar\'e series for subspaces defined by principal ideals. For
plane curve singularities and nondegenerate singularities this Poincar\'e series yields topological and geometric information.
We compare this Poincar\'e series with the one introduced in \cite{sabirebeling}.
In few cases they are equal and we show that the Poincar\'e series we consider in this paper in general yields more information.
}}
\\ \\
${}$
\begin{center}
\textsc{0. Introduction}
\end{center}
${}$
In \cite{cdk} one introduced a Poincar\'e series induced by a filtration on the ring of germs of a complex variety.
This Poincar\'e series has been studied for several kinds of singularities, see for example \cite{triomonodromy}, \cite{trio1}, \cite{trio3}, \cite{ebeling}, \cite{cutkosky}, \cite{lemahieu}, \cite{pedro} and \cite{nemethi}. In some cases this Poincar\'e series determines the topology of the singularity and is related to the zeta function of monodromy of the singularity.

In these works one considers multi-index filtrations defined by valuations on the local ring at the singularity and in \cite{ebeling} and very recently in \cite{sabirebeling} one considers valuations on an ambient smooth space of the singularity that correspond to facets of the Newton polyhedron.
Here we study Poincar\'e series induced by multi-index filtrations coming from arbitrary valuations on the ambient space where at least one of them is centred at the maximal ideal of the local ring of the singularity considered in the ambient space.
In an upcoming paper, we study this Poincar\'e series also for valuations where none of them is centred at the maximal ideal.

The Poincar\'e series we introduce here is defined in an algebraic way and differs from Poincar\'e series studied before in the sense that there is no notion of fibre that corresponds to our Poincar\'e series. We go into this in Section 1. We compute this Poincar\'e series for a subspace corresponding to a principal ideal. A nice A'Campo type formula shows up.
In Section 2 we compare our Poincar\'e series with the one defined in \cite{sabirebeling}. The Poincar\'e series in \cite{sabirebeling} is consistent with a notion of geometric fibre but the price to pay is that this filtration is less richer than the one we consider here. This follows also from the formulas they computed for their Poincar\'e series. When the Newton polytope is bi-stellar, we show that both Poincar\'e series coincide. The property of bi-stellar generalises the notion of stellar introduced in \cite{sabirebeling}.
In Section 3 we study this Poincar\'e series for a plane curve singularity that is a general element in an ideal in $\C\{x,y\}$ where the embedded filtration comes from the Rees valuations of the ideal. Using the result in \cite{trionieuw}, we show that the Poincar\'e series determines and is determined by the embedded topology of the plane curve singularity.
In Section 4 we compare different Poincar\'e series for affine toric varieties.
In Section 5 we study the Poincar\'e series for nondegenerate singularities where the embedded filtration is now induced from the Newton polyhedron. We show that one can recover the Newton polyhedron from the Poincar\'e series and thus, in particular, that the zeta function of monodromy can be deduced from the Poincar\'e series.
\\
${}$
\begin{center}
\textsc{1. Poincar\'e series associated to embedded filtrations}
\end{center}
${}$ \\ Let $(X,o)$ be a germ of a complex analytic space and let
$\mathcal{O}_{X,o}$ be the local ring of germs of functions on
$(X,o)$. Let $\underline{\nu}=\{\nu_1,\cdots,\nu_r\}$ be a set of order functions from $\mathcal{O}_{X,o}$ to
$\mathbb{Z} \cup \{\infty\}$, i.e.\ functions $\nu_j$ that satisfy $\nu_j(f+g) \geq
\textnormal{ min } \{\nu_j(f),\nu_j(g)\}$ and $\nu_j(fg) \geq
\nu_j(f)$, for all $f, g \in \mathcal{O}_{X,o}$, $1 \leq j \leq r$.
An order function $\nu_j$ is called a valuation if moreover it satisfies
$\nu_j(fg) = \nu_j(f)+\nu_j(g)$, for all $f, g \in \mathcal{O}_{X,o}$.
In particular, when $(X,o)$ is
irreducible, $\underline{\nu}$ can be a set of discrete valuations
of the function field $\mathbb{C}(X)$ whose valuation rings contain
$\mathcal{O}_{X,o}$. The set $\underline{\nu}$ defines a multi-index
filtration on $\mathcal{O}_{X,o}$ by the ideals
\[M(\underline{v}):=\{g \in \mathcal{O}_{X,o} \mbox{ $|$ } \nu_j(g)
\geq v_j, 1 \leq j \leq r\}, \qquad \underline{v} \in
\mathbb{Z}^r.\]
If the dimensions of the complex vector spaces $M(\underline{v})/M(\underline{v}+\underline{1})$ are finite for all $\underline{v}
\in \mathbb{Z}^r$, then originally (see \cite{cdk} and \cite{trio1}) the \emph{Poincar\'e series} associated to this multi-index filtration was defined as
\begin{eqnarray*} 
P^{\underline{\nu}}_{X}(t_1,\cdots,t_r):=\frac{\prod_{j=1}^{r}(t_j-1)}{(t_1\cdots t_r - 1)}\sum_{\underline{v}
\in \mathbb{Z}^r} \mbox{dim}(M(\underline{v})/M(\underline{v}+\underline{1}))
\underline{t}^{\underline{v}}.\end{eqnarray*}

Let us now consider an ideal $\I$ in $\mathcal{O}_{X,o}$. We will define a Poincar\'e series associated to
an \emph{embedded filtration} on $\mathcal{O}_{X,o} / \I$. Let $V$ be the analytic subspace of $X$
determined by the ideal $\I$ and let the map $\varphi:
\mathcal{O}_{X,o} \rightarrow \mathcal{O}_{X,o}/\I$ define the
embedding of $V$ in $X$. We set %$J(\underline{v}):=M(\underline{v}) + \I$ and let
$I(\underline{v}):=\varphi(M(\underline{v}) + \I)$. In general the
ideals $I(\underline{v})$ define a multi-index filtration on $\mathcal{O}_{V,o}$ which is
not induced by a set of order functions on
$V$, i.e.\ there do not have to exist order functions $\mu_1,\ldots,\mu_r$
on $V$ such that
$I(\underline{v}) = \{g \in \mathcal{O}_{V,o}  \mid
\mu_j(g) \geq v_j, 1 \leq j \leq r\}$.
The existence of such functions for $\{I(\underline{v})\}$ is equivalent with the condition that
\begin{eqnarray} \label{eqnval}
I(\underline{v_1}) \cap I(\underline{v_2}) = I(\underline{v}),
\end{eqnarray}
where $\underline{v_1}$ and $\underline{v_2}$ are arbitrary tuples in $\Z^r$ and $\underline{v}$ is the tuple of the componentwise maxima of $\underline{v_1}$ and $\underline{v_2}$.
Indeed, if Condition (\ref{eqnval}) holds, then one has order functions $\mu_j$ on $V$ with $\mu_j(g)=v_j$ if and only if $g \in I(0,\ldots,0,v_j,0,\ldots,0) \setminus I(0,\ldots,0,v_j+1,0,\ldots,0)$. 
\begin{definition}\label{def2}
The \emph{Poincar\'e series associated to the embedded multi-index filtration} given by the ideals
$I(\underline{v})$ is the series
\[\mathcal{P}^{\underline{\nu}}_{V}(t_1,\cdots,t_r):=\frac{\prod_{j=1}^{r}(t_j-1)}{(t_1\cdots t_r - 1)}\sum_{\underline{v}
\in \mathbb{Z}^r} \emph{\mbox{dim}}(I(\underline{v})/I(\underline{v}+\underline{1}))
\underline{t}^{\underline{v}}.\]
\end{definition}
\noindent It is useful to notice that $\mbox{dim}(I(\underline{v})/I(\underline{v}+\underline{1})) = \mbox{dim}\left(J(\underline{v})/J(\underline{v}+\underline{1})\right)$, where $J(\underline{v}):=M(\underline{v})+\I$.
If $\{I(\underline{v})\}$ is defined by order functions $\mu_1,\ldots,\mu_r$, then
one has a notion of \emph{fibre} because
to every $g \in \mathcal{O}_{V,o}$ can can attach a value $\underline{\mu}(g)=(\mu_1(g),\ldots,\mu_r(g))$.
Concretely, consider the projection map
\begin{eqnarray*} \label{eqnjv}
j_{\underline{v}}: I(\underline{v}) & \longrightarrow & \nonumber
D_1(\underline{v}) \times \ldots \times D_r(\underline{v}) \\
g & \longmapsto & (a_1(g),\ldots,a_r(g)),
\end{eqnarray*}
where for $1 \leq j \leq r$ the space $D_j(\underline{v})=I(\underline{v})/I(\underline{v}+\underline{e_j})$ where
$\underline{e_j}$ is the $r$-tuple with \emph{j}-th component equal
to $1$ and the other components equal to $0$.
As Condition (\ref{eqnval}) then holds, the following expressions
\begin{eqnarray}
F_{\underline{v}}& := & \textnormal{Im} j_{\underline{v}} \cap
(D_1^{*}(\underline{v}) \times \cdots \times
D_r^{*}(\underline{v})),  \label{eqnfv1}
\\
F_{\underline{v}} &:=
&\left(I(\underline{v})/I(\underline{v+1})\right) \setminus
\bigcup_{i=1}^r(I(\underline{v+e_i})/I(\underline{v+1})) \label{eqnfv2}
\end{eqnarray}
coincide and one has that $\underline{\mu}(g)=\underline{v}$ if and only if
$g + I(\underline{v+1}) \in F_{\underline{v}}$. The space $F_{\underline{v}}$ is invariant with respect to multiplication by nonzero constants; let
$\mathbb{P}F_{\underline{v}}=F_{\underline{v}}/\mathbb{C}^{*}$ be
the projectivisation of $F_{\underline{v}}$.
We then have
\[\mathcal{P}^{\underline{\nu}}_{V}(\underline{t}) = \sum_{\underline{v} \in
\mathbb{Z}^r}\chi(\mathbb{P}F_{\underline{v}})\underline{t}^{\underline{v}},\]
where $\chi(\cdot)$ denotes the Euler characteristic.
If such functions $\mu_1,\ldots,\mu_r$ do not exist for $\{I(\underline{v})\}$, then
(\ref{eqnfv1}) and (\ref{eqnfv2}) are possible generalisations of the notion of fibre
but we then lose some geometric meaning.

In this paper we will investigate the Poincar\'e
series as in Definition \ref{def2} and we will see
in Example \ref{anders} that this Poincar\'e series is different from the two possible
ones induced by the `fibres' (\ref{eqnfv1}) and (\ref{eqnfv2}). \\
\indent We now compute this Poincar\'e series for $\I=(h)$ a
principal ideal.
\begin{theorem} \label{thmembedded}
Let $(X,o)$ be irreducible, $\I=(h)$ a principal ideal in
$\mathcal{O}_{X,o}$ and $V$ the analytic subspace of $(X,o)$ determined by the ideal $\I$. Let $\underline{\nu}=\{\nu_1,\ldots,\nu_r\}$ be a set of discrete valuations of
$\mathbb{C}(X)$ centred at the maximal ideal of
$\mathcal{O}_{X,o}$. We write
$\underline{q}=\underline{\nu}(h)$. Then
\begin{eqnarray}\label{compu}
\mathcal{\mathcal{P}}^{\underline{\nu}}_{V}(\underline{t})=(1-\underline{t}^{\underline{q}})P^{\underline{\nu}}_X(\underline{t}).
\end{eqnarray}
\end{theorem}
\noindent \emph{Proof.} \quad  For a set $A
\subset \{1,\cdots,r\}$, let $\underline{1}_A$ be the r-tuple with
j-th component equal to 1 if $j \in A$ and else equal to 0, $1 \leq
j \leq r$.
%\begin{eqnarray*}
%\alpha_A:  \mathbb{Z}^s & \longrightarrow &\mathbb{Z}^s \\
%\underline{v} & \longmapsto & \underline{v}'
%\end{eqnarray*}
%where $v'_i=v_i-1$ if $i \in A$ and $v'_i=v_i$ if $i \notin A$.
The coefficient of $\underline{t}^{\underline{v}}$ in the left hand
side of (\ref{compu}) is
\begin{eqnarray*}
 (-1)^{r+1} \sum_{A \subset \{1,\cdots,r\}} (-1)^{\#A} \mbox{dim}
\frac{\mathcal{O}_{X,o}}{J(\underline{v}-\underline{1}_A+
\underline{1})}.
\end{eqnarray*}
The coefficient of $\underline{t}^{\underline{v}}$ in the right hand
side of (\ref{compu}) is equal to
\begin{eqnarray*}
 (-1)^{r+1} \sum_{A \subset \{1,\cdots,r\}} (-1)^{\#A} \mbox{dim}
\frac{\mathcal{O}_{X,o}}{M(\underline{v}-\underline{1}_A+
\underline{1})}-\\(-1)^{r+1} \sum_{A \subset \{1,\cdots,r\}}
(-1)^{\#A} \mbox{dim}
\frac{\mathcal{O}_{X,o}}{M(\underline{v}-\underline{1}_A-\underline{q}+
\underline{1})}.
\end{eqnarray*}
These sums can be rewritten such that for $A \subset \{1,\cdots,r\}$, it is enough to prove for some $j \in A$
that \begin{eqnarray*}& & \textnormal{dim
}\frac{J(\underline{v}-\underline{1}_A+
\underline{1})}{J(\underline{v}-\underline{1}_{A\setminus
\{j\}}+ \underline{1})} \\ & =  & \textnormal{dim
}\frac{M(\underline{v}-\underline{1}_A+
\underline{1})}{M(\underline{v}-\underline{1}_{A\setminus
\{j\}}+ \underline{1})}- \textnormal{dim
}\frac{M(\underline{v}-\underline{1}_A-\underline{q}+
\underline{1})}{M(\underline{v}-\underline{1}_{A\setminus
\{j\}}-\underline{q}+ \underline{1})}.\end{eqnarray*}
This follows immediately from the fact that the kernel of the
projection map
\[\frac{M(\underline{v}-\underline{1}_{A}+
\underline{1})}{M(\underline{v}-\underline{1}_{A\setminus
\{j\}}+ \underline{1})} \longrightarrow
\frac{J(\underline{v}-\underline{1}_{A}+
\underline{1})}{J(\underline{v}-\underline{1}_{A\setminus
\{j\}}+ \underline{1})}\] is equal to
$$\frac{(h)M(\underline{v}-\underline{1}_{A}-\underline{q}+
\underline{1})}{(h)M(\underline{v}-\underline{1}_{A\setminus
\{j\}}-\underline{q}+ \underline{1})}.$$ \hfill
$\blacksquare$
${}$
\begin{example} \label{anders}
\emph{Let $X=\C^2$ and $h(x,y)=x^6y^2+y^8$. We choose monomial
valuations $\nu_1$ and $\nu_2$ on $\mathcal{O}_{X,o}$ given by $\nu_1(x^ay^b)=2a+3b$ and
$\nu_2(x^ay^b)=4a+3b$. One can compute that the coefficient of
$t_1^{20}t_2^{28}$ in
$(1-\underline{t}^{\underline{q}})P^{\underline{\nu}}_X(\underline{t})$
is $0$ although both fibres $F_{(20,28)}$ as defined in (\ref{eqnfv1}) and (\ref{eqnfv2}) contain
the monomial $x^4 y^4$. \hfill $\square$}
\end{example}
${}$ \\ ${}$
\vspace*{-0.5cm}
\begin{center}
\textsc{2. Geometric embedded filtrations versus algebraic embedded filtrations.}
\end{center}
${}$
\\
Let $h: (\C^d,o) \rightarrow (\C,o)$ be a germ of a holomorphic function. Recently in \cite{sabirebeling} one considered a Poincar\'e series on $\mathcal{O}_{V,o}=\mathcal{O}_{\C^d,o}/(h)$, induced by a \emph{Newton filtration}. We write $h(\underline{x}) = \sum_{\underline{k}
\in \Z_{\geq 0}^d} a_{\underline{k}} \underline{x}^{\underline{k}},$ where $\underline{k}=(k_1, \ldots, k_d)$ and $\underline{x}^{\underline{k}} =
x_1^{k_1} \cdot \ldots \cdot x_d^{k_d}.$ The \emph{support of $h$} is
$\textnormal{supp}\, h := \{\underline{k} \in \Z_{\geq 0}^d \, | \, a_{\underline{k}} \neq 0
\}.$ The \emph{Newton polyhedron of $h$ at the origin} is
the convex hull in $\mathbb{R}_{\geq 0}^d$ of $\bigcup_{\underline{k} \in
\textnormal{{supp}}\, h} \underline{k} + \mathbb{R}_{\geq 0}^d$ and the \emph{Newton polytope of $h$ at the origin}
is the compact boundary of the Newton polyhedron of $h$ at the origin.

Let $\nu_1,\ldots,\nu_r$ be the monomial valuations on
$\mathcal{O}_{\C^d,o}$ corresponding to the facets of the Newton polytope of $h$, i.e.\ for a compact facet $\tau$ with the affine space through $\tau$ given by the equation $a_1x_1+\cdots + a_dx_d=N_{\tau}$, the corresponding valuation $\nu$ acts as follows: $\nu(x_1^{m_1}\ldots x_d^{m_d})=a_1m_1+\cdots + a_dm_d$.
W. Ebeling and S. M. Gusein-Zade study the Poincar\'e series $P_{\{\omega_i\}}(\underline{t})$ on $V$ induced by the
order functions
\[\omega_i(g):= \textnormal{max }\{\nu_i(g') \mid g'-g \in (h)\}.\]
On the other hand we can consider the Poincar\'e series as defined in Definition \ref{def2}, with $\I=(h)$.
The Poincar\'e series in \cite{sabirebeling} has a more geometric meaning because one has fibres $F_{\underline{v}}$ and
$P_{\{\omega_i\}}(\underline{t})=\sum_{\underline{v} \in
\mathbb{Z}^r}\chi(\mathbb{P}F_{\underline{v}})\underline{t}^{\underline{v}}$.
However, our Poincar\'e series - which is rather algebraic - contains more information about the singularity $(V,o)$. In general, the embedded filtration we introduce in this article is richer because $I(\underline{v})$ is not necessarily determined by $I(v_1,0,\ldots,0), I(0,v_2,0$ $,\ldots,0), \ldots, I(0,\ldots,0,v_d)$, whereas the ideals that appear in the Poincar\'e series $P_{\{\omega_i\}}(\underline{t})$ are.
We comment further on this in this section.
\\
In \cite{sabirebeling} a Newton polytope is called \emph{stellar} if all its facets have a common vertex.
If the Newton polytope of $h$ is stellar, then they show (\cite[Theorem 2]{sabirebeling}) that
\[P_{\{\omega_i\}}(\underline{t})=(1-\underline{t}^{\underline{\nu}(h)})P^{\underline{\nu}}_X(\underline{t}).\]
However, if the Newton polytope is not stellar, the information on $(V,o)$ can be lost in the Poincar\'e series $P_{\{\omega_i\}}(\underline{t})$. Indeed, for $h$ a germ of a holomorphic function on $(\C^2,o)$ with $o$ an isolated critical point of $h$, they prove (\cite[Theorem 1]{sabirebeling}) that
\[P_{\{\omega_i\}}(\underline{t})=P^{\underline{\nu}}_X(\underline{t}).\]
\begin{proposition}
If the ideals $I(\underline{v})$ satisfy Condition (\ref{eqnval}), then
$P_{\{\omega_i\}}(\underline{t})=\mathcal{P}^{\underline{\nu}}_{V}(\underline{t})$.
\end{proposition}
\noindent \emph{Proof.} \quad
%Notice that the ideals $J(\underline{v})$ satisfy Condition (\ref{eqnval}) if and only if
%the ideals $I(\underline{v})$ satisfy Condition (\ref{eqnval}).
It is sufficient to verify that then
\[I(0,\ldots,0,v_i,0,\ldots,0) = \{g \in \mathcal{O}_{V,o} \mid \omega_i(g) \geq v_i\}.\vspace*{-0.5cm} \]
\hfill $\blacksquare$   
${}$
\\ \\
We will now characterise the germs $h$ of holomorphic functions on $\C^d$ for which Condition (\ref{eqnval}) is satisfied, and thus for which these geometric and algebraic Poincar\'e series coincide.
Obviously this will depend on the Newton polytope of $h$.
We will use the following lemma to give the characterisation.
\begin{lemma} \label{lemmatechnical}
Let $\{M(\underline{v})\}$ be ideals in a local ring $R$ that satisfy Condition (\ref{eqnval}) and consider an ideal $\I$ in $R$. Then for $\underline{v_1}, \underline{v_2} \in \Z^r$, the following conditions are equivalent:
\begin{enumerate}
\item $(M(\underline{v_1})+\I) \cap (M(\underline{v_2})+\I) = M(\underline{v})+\I$ where  $\underline{v}$ is the tuple of the componentwise maxima of $\underline{v_1}$ and $\underline{v_2}$;
\item $(M(\underline{v_1}) + M(\underline{v_2}))\cap \I = (M(\underline{v_1}) \cap \I) + (M(\underline{v_2}) \cap \I)$.
\end{enumerate}
\end{lemma}
\noindent \emph{Proof.} \quad
Suppose that the first equation of sets holds. Take $f \in (M(\underline{v_1}) + M(\underline{v_2})) \cap \I$. Then  $f=f_1+f_2$ with $f \in \I$, $f_1 \in M(\underline{v_1})$ and $f_2 \in M(\underline{v_2})$. Thus $f_1 \in (M(\underline{v_1}) + \I) \cap (M(\underline{v_2}) + \I)$.
By the hypothesis it then follows that $f_1 \in M(\underline{v})+\I$ and so we can write
$f_1= g + k$ with $g \in M(\underline{v})$ and $k \in \I$. We have $f = (f_1-g)
+ (f_2+g)$ and $f_1-g \in M(\underline{v_1}) \cap \I$ and $f_2+g \in M(\underline{v_2}) \cap \I$.
Hence $f \in (M(\underline{v_1}) \cap \I)+(M(\underline{v_2}) \cap \I)$.
\\ \\
We now suppose that the second equation of sets holds. %Let
%\[Q_\I(\underline{w_1},\underline{w_2}):=[(M(\underline{w_1}) + M(\underline{w_2}))\cap \I]/ [(M(\underline{w_1}) \cap \I) + (M(\underline{w_2}) \cap \I)].\]
We consider the following exact sequences:
\footnotesize{ \begin{eqnarray}
 0 \rightarrow  (M(\underline{v_1}) \cap \I)/(M(\underline{v}) \cap \I) \rightarrow M(\underline{v_1})/M(\underline{v}) \rightarrow (M(\underline{v_1})+\I)/(M(\underline{v})+\I) \rightarrow  0 \label{exact1}\\
0  \rightarrow  (M(\underline{v_2}) \cap \I)/(M(\underline{v}) \cap \I) \rightarrow M(\underline{v_2})/M(\underline{v}) \rightarrow (M(\underline{v_2})+\I)/(M(\underline{v})+\I) \rightarrow 0 \label{exact2}%\\
%0 & \rightarrow & (M(\underline{w_1}) \cap \I)/ (M(\underline{w}) \cap \I) \rightarrow
%[(M(\underline{w_1}) + M(\underline{w_2}))\cap \I]/ (M(\underline{w_2}) \cap \I) \rightarrow  \nonumber \\ & & Q_\I(\underline{w_1},\underline{w_2}) \rightarrow 0  \label{exact3}\\
%0 & \rightarrow & (M(\underline{w_2}) \cap \I)/ (M(\underline{w}) \cap \I) \rightarrow
%[(M(\underline{w_1}) + M(\underline{w_2}))\cap \I]/ (M(\underline{w_1}) \cap \I) \rightarrow \nonumber \\ & & Q_\I(\underline{w_1},\underline{w_2}) \rightarrow 0 \label{exact4}
\end{eqnarray}}
\normalsize
We take a set $B'_i$ of elements in $M(\underline{v_i}) \cap \I$ that give rise to a basis of the vector space
$(M(\underline{v_i}) \cap \I)/(M(\underline{v}) \cap \I)$ , $i \in \{1,2\}$.
For $i \in \{1,2\}$, we use the exact sequences (\ref{exact1}) and (\ref{exact2}) to add to $B'_i$ a set of elements $B''_i$ of $M(\underline{v_i})$ whose images are a basis of the quotient
$(M(\underline{v_i})+\I)/(M(\underline{v})+\I)$. So the classes of the elements of
$B_i:= B'_i \cup B''_i$ form a basis of $M(\underline{v_i})/M(\underline{v})$, $i=1,2$.
We also have that the classes of the elements in
$B_1 \cup B_2$ form a system of generators for the vector space $(M(\underline{v_1}) + M(\underline{v_2}))/
M(\underline{v})$. Analogously, the classes of the elements in
$B'_1 \cup B'_2$ are a system of generators for the vector space $ [(M(\underline{v_1}) \cap \I) + (M(\underline{v_2}) \cap \I)]/(M(\underline{v}) \cap \I)$.

Condition (\ref{eqnval}) implies that the set $B_1 \cup B_2$ consists of elements
of $M(\underline{v_1}) + M(\underline{v_2})$ whose classes module $M(\underline{v})$ are linearly independent and hence the classes of the elements in $B_1 \cup B_2$ are a basis of $(M(\underline{v_1}) + M(\underline{v_2}))/
M(\underline{v})$. Indeed, if $x$ is a linear combination of elements in $B_1$ and $y$ is a linear combination of elements of $B_2$ and if $x+y=z$ with $z \in M(\underline{v})$, then
$x=z-y \in M(\underline{v_1}) \cap M(\underline{v_2})$. Condition (\ref{eqnval}) then implies that $x, y \in M(\underline{v})$.

%The hypothesis (\ref{eqnval}) implies that $Q_\I(\underline{w_1},\underline{w_2})=0$. The exact sequences (\ref{exact3}) and (\ref{exact4}) allow us to deduce that the set $B'_1$ (resp. $B'_2$) is also a set of elements in
%$(M(\underline{w_1}) + M(\underline{w_2}))\cap \I$ that gives rise to a basis for the quotient
%$[(M(\underline{w_1}) + M(\underline{w_2}))\cap \I]/(M(\underline{w_2}) \cap \I)$ (resp. $[(M(\underline{w_1}) + M(\underline{w_2}))\cap \I]/(M(\underline{w_1}) \cap \I)$).
Let us now take an element $f \in (M(\underline{v_1})+\I) \cap (M(\underline{v_2})+\I)$, so $f= f_1+k_1=f_2+k_2$, $f_i \in M(\underline{v_i})$ and $k_i \in \I$, $i=1,2$. For $i=1,2$, we can write $f_i = f'_i+f''_i
+ m_i$ with $f'_i$ a linear combination of elements in $B'_i$, $f''_i$ a linear combination of elements in $B''_i$ and $m_i \in M(\underline{v})$. Then $g := k_2–-k_1= f_1-f_2 \in (M(\underline{v_1}) + M(\underline{v_2}))\cap \I = (M(\underline{v_1}) \cap \I) + (M(\underline{v_2}) \cap \I)$ by the hypothesis. Hence $g= g'_1 + g'_2 + m$, with $g'_i$
a linear combination of elements of $B'_i$, $i=1,2$, and $m \in M(\underline{v}) \cap \I$. We obtain
\[(g'_1-f'_1) –- f''_1+(g'_2+f'_2) +f''_2= m_1-m_2-m.\]
The right hand side of this equality is contained in $M(\underline{v})$. The left hand side is a sum of four terms that are linear combinations of elements respectively from $B'_1, B''_1,B'_2$ and $B''_2$.
As the classes of the elements in $B_1 \cup B_2$ are a basis of $(M(\underline{v_1}) + M(\underline{v_2}))/
M(\underline{v})$, we get in particular $f''_1=0$. Hence $f= m_1+(f'_1+k_1) \in M(\underline{v})+\I$.
\hfill $\blacksquare$
\\${}$
\begin{definition}
A Newton polytope is called \emph{bi-stellar} if every two facets of the Newton polytope have a non-empty intersection.
\end{definition}
\begin{proposition}
The Newton polytope of $h$ is bi-stellar if and only if the ideals $M(\underline{v})+(h)$ satisfy Condition (\ref{eqnval}).
\end{proposition}
\noindent \emph{Proof.} \quad
Say the Newton polytope of $h$ has $r$ compact facets inducing the monomial valuations $\nu_1,\ldots,\nu_r$ on $\C^d$.
Suppose that the Newton polytope of $h$ is bi-stellar. %This is equivalent to say that there exists a monomial $\underline{x}^{\underline{a}}:=x_1^{a_1}\cdots x_n^{a_n}$ in the support of $h$ such that $\underline{\nu}(h)=\underline{\nu}(\underline{x}^{\underline{a}})$. We denote the value $\underline{\nu}(h)$ by $\underline{q}$.
By Lemma \ref{lemmatechnical}, it suffices to show that for all $\underline{v_1}, \underline{v_2} \in \Z^r$ one has that
\[(I(\underline{v_1}) \cap (h)) + (I(\underline{v_2}) \cap (h)) = (I(\underline{v_1})+I(\underline{v_2}))\cap (h).\]
Let $gh \in I(\underline{v_1})+I(\underline{v_2}))$ and $\underline{q}:=\underline{\nu}(h)$. We write
$g=g_1+g_2$ with $g_1=\sum \lambda_{\underline{a}}\underline{x}^{\underline{a}}$ and
$\underline{x}^{\underline{a}} \notin I(\underline{v_1}-\underline{q})+I(\underline{v_2}-\underline{q})$, for all
$\underline{x}^{\underline{a}}$ in supp $g_1$, and
$g_2=\sum \lambda_{\underline{b}}\underline{x}^{\underline{b}}$ with
$\underline{x}^{\underline{b}} \in I(\underline{v_1}-\underline{q})+I(\underline{v_2}-\underline{q})$, for all
$\underline{x}^{\underline{b}}$ in supp $g_2$.
Suppose that $g_1 \neq 0$. We take a monomial $\underline{x}^{\underline{a}}$ in supp $g_1$. Then there exist $i,j \in \{1,\ldots,r\}$ such that $\nu_i(\underline{x}^{\underline{a}}) < v_{1,i}-q_i$ and $\nu_j(\underline{x}^{\underline{a}}) < v_{2,j}-q_j$.

We first consider the case where $i \neq j$.
Let $N$ be the set of the monomials in supp $g_1$ that are minimal for the pair $(\nu_i,\nu_j)$, i.e.\ $\underline{x}^{\underline{c}} \in N$ if and only if there does not exist a monomial $\underline{x}^{\underline{d}}$ in supp $g_1$ for which $\nu_i(\underline{x}^{\underline{d}})< \nu_i(\underline{x}^{\underline{c}})$ and
$\nu_j(\underline{x}^{\underline{d}})< \nu_j(\underline{x}^{\underline{c}})$.
Let $M$ be the set of monomials $\underline{x}^{\underline{m}}$ in supp $h$ for which $\nu_i(\underline{x}^{\underline{m}})=q_i$ and $\nu_j(\underline{x}^{\underline{m}})=q_j$. As the Newton polytope is bi-stellar, $M$ is not empty. For the monomials $\underline{x}^{\underline{m}}\underline{x}^{\underline{c}}$ with
$\underline{x}^{\underline{m}} \in M$ and $\underline{x}^{\underline{c}} \in N$, we thus have that $\nu_i(\underline{x}^{\underline{m}}\underline{x}^{\underline{c}}) < v_{1,i}$ and $\nu_j(\underline{x}^{\underline{m}}\underline{x}^{\underline{c}}) < v_{2,j}$. As such monomials do not appear in the support of $gh$ they should be canceled. It follows that at least one such monomial $\underline{x}^{\underline{m}}\underline{x}^{\underline{c}}$ has to be equal to a monomial $\underline{x}^{\underline{h}}\underline{x}^{\underline{a'}}$ with $\underline{x}^{\underline{h}} \in \textnormal{supp } h \setminus M$ and $\underline{x}^{\underline{a'}}$ in supp $g_1$. Say $\nu_i(\underline{x}^{\underline{h}}) > q_i$. We then find that $\nu_i(\underline{x}^{\underline{a'}}) < \nu_i(\underline{x}^{\underline{c}})$ and $\nu_j(\underline{x}^{\underline{a'}}) \leq \nu_j(\underline{x}^{\underline{c}})$ which contradicts the fact that  $\underline{x}^{\underline{c}} \in N$ and thus $g_1=0$.

Suppose now that $i=j$. Let $N$ be the set of monomials with support in $g_1$ that are minimal for the valuation $\nu_i$, i.e.\ $\underline{x}^{\underline{c}} \in N$ if and only if there does not exist a monomial $\underline{x}^{\underline{d}}$ in supp $g_1$ for which $\nu_i(\underline{x}^{\underline{d}})< \nu_i(\underline{x}^{\underline{c}})$. Let $M$ be the set of monomials $\underline{x}^{\underline{m}}$ in the support of $h$ for which $\nu_i(\underline{x}^{\underline{m}})=q_i$. Analogously there then has to be a monomial
$\underline{x}^{\underline{m}}\underline{x}^{\underline{c}}$ with $\underline{x}^{\underline{m}} \in M$ and $\underline{x}^{\underline{c}} \in N$ that is equal to a monomial $\underline{x}^{\underline{h}}\underline{x}^{\underline{a'}}$, with
$\underline{x}^{\underline{h}} \in \textnormal{supp } h \setminus M$ and $\underline{x}^{\underline{a'}} \in$ supp $g_1$.
Again we get a contradiction because $\underline{x}^{\underline{c}}$ would not be minimal for $\nu_i$.
%We denote $h'$ for the part of $h$ in which every monomial has value $\underline{q}$, so $\underline{\nu}(h')=\underline{q}$ and $\nu_j(h-h') > q_j$, for every $j \in \{1,\ldots,r\}$.
%Then $g_1h'$ can not be canceled in $gh$. This gives a contradiction to the fact that $gh \in I(\underline{w_1})+I(\underline{w_2})$.
\\ \\
We now suppose that the ideals $M(\underline{v})+(h)$ satisfy Condition (\ref{eqnval}).
If the Newton polytope of $h$ would not be bi-stellar, then there would exist two valuations $\nu_i,\nu_j$ with $i,j\in \{1,\ldots,r\}$ for which the sets of monomials $M_i=\{\underline{x}^{\underline{m}} \in \textnormal{supp}(h) \mid \nu_i(\underline{x}^{\underline{m}})=q_i\}$ and $M_j=\{\underline{x}^{\underline{m}} \in \textnormal{supp}(h) \mid \nu_j(\underline{x}^{\underline{m}})=q_j\}$ would be disjoint. Let $h_i$ be the part of $h$ with support in $M_i$, so $\nu_i(h_i)=q_i$ and $\nu_i(h-h_i)>q_i$. Then $h=h_i+(h-h_i)$ with $\nu_i(h_i) < \nu_i(h-h_i)$,
$\nu_j(h-h_i) < \nu_j(h_i)$ and $\underline{\nu}(h) \leq \underline{\nu}(h_i)$ and $\underline{\nu}(h) \leq \underline{\nu}(h-h_i)$. By Lemma \ref{lemmatechnical} it follows that $h$ would be contained
in $(M(\underline{\nu}(h_i))+M(\underline{\nu}(h-h_i))) \cap (h)$ but not in $(M(\underline{\nu}(h_i)) \cap (h)) + (M(\underline{\nu}(h-h_i))) \cap (h))$, contradicting Condition (\ref{eqnval}).
\hfill $\blacksquare$
\\
${}$
\begin{center}
\textsc{3. Embedded filtrations for plane curve singularities}
\end{center}
${}$
\\
Let $X=\mathbb{C}^2$ and let $\I$ be a primary
ideal in $\mathcal{O}_{X,o}$. If $\phi: Z \rightarrow \C^2$
is a principilisation of $\I$, then $\phi$ is realised by blowing up a constellation of points
$\{Q_\sigma\}_{\sigma \in G}$. The map $\phi$ factorises through
the normalised blowing up of $\I$ which we will denote by
$\overline{Bl_\I(\C^2)}$. Let $\sigma$ be the morphism $Z \rightarrow
\overline{Bl_{\I}(\C^2)}$ in this factorisation. For $\sigma \in G$, we denote the exceptional divisor of the blowing-up in $Q_\sigma$ by $E_\sigma$, as well as its strict transform under following blowing-ups, and $\mathcal{D}:=\cup_{\sigma \in G}E_{\sigma}$.
Blowing up a point $Q_\sigma$ induces a discrete valuation $\nu_\sigma$ on $\C(X) \setminus \{0\}$: for $g \in \C(X) \setminus \{0\}$, the value $\nu_\sigma(g)$ is the order of the pullback of g along $E_\sigma$.
The valuation $\nu_\sigma$ is called Rees for $\I$ if its
centre in $\overline{Bl_{\I}(X)}$ is a divisor. We have that $\nu_\sigma$ is Rees for $\I$ if and only if
the strict transform of a general element in $\I$ intersects $E_\sigma$ (see for example \cite[Lemma 8]{holomorfie}).
Say that $E_1,\ldots,E_r$ are the exceptional components that give rise to Rees valuations
$\nu_1,\ldots,\nu_r$.

Let us now consider a general element $h$ in the ideal $\I$ and let $V$ be the hypersurface given by $\{h=0\}$.
In \cite{trio1} one studied the Poincar\'e series $P^{\underline{\nu}}_V(\underline{t})$ that is defined by the filtration on $\mathcal{O}_{V,o}$ induced by the essential valuations $\underline{\nu}$ of the minimal resolution of the plane curve $V$. One showed that that Poincar\'e series
contained the same information as the embedded topology of the curve and that $P^{\underline{\nu}}_V(t,\ldots,t)$ equals the zeta function of monodromy.

We will now study the Poincar\'e series of the embedded filtration on $\mathcal{O}_{V,o}$ induced by the Rees valuations $\underline{\nu}=(\nu_1,\ldots,\nu_r)$.
For $1 \leq j \leq r$, suppose that $E_j$ is intersected $n_j$ times by
the strict transform of $\{h=0\}$. Let $E^{\bullet}_{\sigma}$ be
$E_{\sigma}$ without the intersection points with the other
components of $\mathcal{D}$ and let $E^{\circ}_{\sigma}$ be
$E_{\sigma}$ without the intersection points of the other components
of $\phi^{-1}(h^{-1}\{o\})$. Let $I$ be the intersection matrix of
the $\{E_{\sigma}\}_{\sigma \in G}$ and let $M=-I^{-1}$. Let $\mathcal{C}_{\sigma}$ be
a curvette through $E_{\sigma}$ (i.e.\ the projection by $\phi$ of a
smooth curve transversal to $E_{\sigma}$ and not intersecting other
components of $\mathcal{D}$). The entry $m_{\sigma,\tau}$ in $M$ is
then also equal to $\nu_{\tau}(\mathcal{C}_{\sigma})$. %Let us denote the total
%transform of $\mathcal{C}_{\sigma}$ by
%$\tilde{\mathcal{C}_{\sigma}}$ and its strict transform by
%$\mathcal{C}_{\sigma}'$. If $D_j$ is the Cartier divisor such that
%$\tilde{\mathcal{C}_{\sigma}}=\mathcal{C}_{\sigma}'+D_j$, then
%$D_j=\sum_{i=1}^r \nu_{i}(\mathcal{C}_{j}) E_i$. Moreover
%$D_i.E_j=-\delta_{i,j}$ such that
%$\nu_{i}(\mathcal{C}_{j})=-D_i.D_j$. We have
%$\mathcal{C}_{i}.\mathcal{C}_{j}=(\mathcal{C}_{i}+D_i).(\mathcal{C}_{j}+D_j)=
%D_i.D_j$ as $\mathcal{C}_{i}$ and $\mathcal{C}_{j}$ are curvettes.
%If $e_{i,k}$ is the multiplicity of $\mathcal{C}_{i}$ in the point
%$Q_k$, then $m_{i,j}=\sum_k e_{i,k}e_{j,k}$. \\
%We compute the
%Poincar\'e series induced by the multi-index filtration by the
%ideals
%\[J(\underline{v}):=\{ g \in \mathcal{O}_{X,o} \mbox{ $|$ } \nu_j(g) \geq v_j, 1 \leq j \leq r\} + (h), \quad \underline{v} \in \mathbb{Z}^r.\]
\begin{theorem} \label{thmcurvas}
The Poincar\'e series $\mathcal{P}^{\underline{\nu}}_V(\underline{t})$ determines and is determined by the embedded topology of $\{h=0\}$.
\end{theorem}
\noindent \emph{Proof.} \quad
By Theorem \ref{thmembedded},
$\mathcal{P}^{\underline{\nu}}_V(t_1,\cdots,t_r)=(1-t_1^{q_1}\cdots
t_r^{q_r})P^{\underline{\nu}}_X(t_1,\cdots,t_r)$, with $\underline{q}=\underline{\nu}(h)$. The Poincar\'e
series $P^{\underline{\nu}}_X(\underline{t})$ induced by plane divisorial valuations is computed in
\cite{duo}. For $\sigma \in G$, let $\underline{m_{\sigma}}=(m_{\sigma,1},\ldots,m_{\sigma,r})$. Then one has
\[P^{\underline{\nu}}_X(t_1,\ldots,t_r)=\prod_{\sigma \in G} (1-\underline{t}^{\underline{m}^{\sigma}})^{-\chi(E^{\bullet}_{\sigma})}.\]
As $E_j$ is intersected $n_j$ times by
the strict transform of $\{h=0\}$, $1 \leq j \leq r$, we get
$q_i=\sum_{j=1}^r n_j m_{i,j}$, for $1 \leq i \leq r$. If the curve
is irreducible (i.e.\ $r=1$) and if $E_1$ is intersected by the strict
transform then $n_1=1$ and $q_1=m_{1,1}$. We then have
\[\mathcal{P}^{\underline{\nu}}_V(t)=(1-t^{m_{1,1}})\prod_{\sigma \in
G}(1-t^{m_{\sigma,1}})^{-\chi(E^{\bullet}_{\sigma})}=\prod_{\sigma
\in G}(1-t^{q_{\sigma}})^{-\chi(E^{\circ}_{\sigma})}.\] Hence it
follows that $\mathcal{P}^{\underline{\nu}}_V(t)$ is then equal to the zeta function of monodromy $\zeta_V(t)$ (\cite{A'Campo}).

Suppose now that the curve is reducible. We will show that the factor
$(1-\underline{t}^{\underline{q}})$ can not be canceled by a factor
$(1-\underline{t}^{\underline{m}^{\sigma}})^{-\chi(E^{\bullet}_{\sigma})}$
of $P^{\underline{\nu}}_X(\underline{t})$. As $q_i=\sum_{j=1}^r n_j m_{i,j}$,
we have that $q_i > m_{i,j}$ for all $j \in \{1,\cdots,r\}$. If $\nu_{\sigma}$ is not a Rees valuation,
then there exists always a valuation $\nu_j$ which is
Rees and such that $Q_j$ lies above $Q_{\sigma}$. Then $m_{\sigma,j} < m_{j,j}$ %(indeed, curvettes satisfy the
%linear proximity equalities, such that $e_{i,k} \leq e_{j,k}$ when
%$j \rightarrow i$)
and thus $q_i > m_{\sigma,i}$ for all $\sigma \in G$.
Thus $\underline{q}$ is the biggest exponent in the cyclotomic factors in $\mathcal{P}^{\underline{\nu}}_V(\underline{t})$.
This makes that we can extract the value $\underline{q}$ and the Poincar\'e series $P^{\underline{\nu}}_X(\underline{t})$ from the Poincar\'e series $\mathcal{P}^{\underline{\nu}}_V(\underline{t})$.
In \cite{trionieuw} it has been shown that $P^{\underline{\nu}}_X(\underline{t})$ determines the
dual graph of the divisors $\{E_{\sigma}\}_{\sigma \in G}$ and thus the matrix $M$. As $M$ is invertible, it follows that we can now compute the
numbers $n_j$, $1 \leq j \leq r$. Hence the dual resolution graph of $\{h=0\}$ is known and
the Poincar\'e series $\mathcal{P}^{\underline{\nu}}_V(\underline{t})$ determines the embedded topology of $\{h=0\}$. \hfill $\blacksquare$
\\ \\ For the Poincar\'e series we study here we get that the zeta function of monodromy
\[\zeta_V(t)=\frac{\mathcal{P}^{\underline{\nu}}_V(t^{n_1},\ldots,t^{n_r})\prod_{j=1}^r (1-t^{q_j})^{n_j}}{(1-t^{\sum_{j=1}^r n_j q_j})}.\]

\begin{remark}
\emph{If $V=\{h=0\}$ is a reduced plane curve singularity and $Z \rightarrow \C^2$ is a
concrete embedded resolution of singularities for $V$, then the function $h$ becomes a
general element for some convenient primary ideal $\I$ such that $\I\mathcal{O}_Z$ is locally principal.
Hence Theorem \ref{thmcurvas} can be applied to any reduced curve singularity and a chosen embedded
resolution for it and, in particular, for the minimal one.}
\end{remark}

%\begin{corollary}
%Let $I$ be an ideal in $\mathbb{C}[x,y]$ and let $C$ be a generic
%curve in $I$. The embedded Poincar\'{e} series for $C$ with respect
%to the Rees valuations of $I$ determines the zeta function of
%monodromy and the topological zeta function.
%\end{corollary}
%\emph{Proof.} \quad This is a direct consequence of the fact that
%one can deduce the embedded resolution graph out of $P$ and then
%follows by A'Campo's formula and by definition of the topological
%zeta function. \hfill $\blacksquare$
${}$
\newpage
\begin{center}
\textsc{4. Poincar\'e series for toric varieties}
\end{center}
${}$ \\
%Let $\underline{\mu}=(\mu_1,\ldots,\mu_r)$ be monomial valuations on $\C(x_1,\ldots,x_d)$ with centre in the maximal ideal.
%If $\I=(f_1,\ldots,f_p) \subset \mathcal{O}_{\C^d,o}$ is a monomial ideal or a quasi-homogeneous ideal for the weights
%$\mu_1,\ldots,\mu_r$, then it is easy to verify that the ideals $M(\underline{v}) +\I$ satisfy the condition ($2$)
%of Lemma \ref{lemmatechnical}.
%If $\I$ is a complete intersection, $\underline{q_i}=\underline{\mu}(f_i)$ ($1 \leq i \leq p$) and $V=$ Spec $(\mathcal{O}_{\C^d,o}/\I)$, then
%\[\mathcal{P}^{\underline{\mu}}_{V}(\underline{t})=(1-\underline{t}^{\underline{q_1}})\ldots (1-\underline{t}^{\underline{q_p}})P^{\underline{\mu}}_{\C^d}(\underline{t}).\]
We now consider the particular case where $V$ is an affine
toric variety that is a complete intersection.
Let $S$ be a semigroup in $M
\cong \mathbb{Z}^d$ such that $S+(-S)=M$ and $S \cap (-S)= {0}$ and let $V=$ Spec
$\mathbb{C}[S]$ be the associated affine toric variety. Let
$\check{\sigma}$ be the cone generated by $S$ and let
$\{s_1,\cdots,s_{d+p}\}$ be a system of generators of
$S=\check{\sigma} \cap M$. Suppose that the embedding of $V$ in
$\mathbb{C}^{d+p}$ is given by the map
\begin{eqnarray*}
\varepsilon: \mathbb{C}[x_1,\cdots,x_{d+p}] & \rightarrow &
\mathbb{C}[S] \\
x_k & \mapsto & \chi^{s_k}.
\end{eqnarray*}
We set $\mbox{deg}(x_k)=s_k$, $1 \leq k \leq d+p$. The toric complete intersection
$V$ is given by an ideal generated by binomials
$h^i=\underline{x}^{\underline{\alpha^i}}-\underline{x}^{\underline{\beta^i}}
\in \mathbb{C}[x_1,\cdots,x_{d+p}]$, $1 \leq i \leq p$. Moreover
deg$(\underline{x}^{\underline{\alpha^i}})=$ deg$(\underline{x}^{\underline{\beta^i}})$
and supp$(\underline{x}^{\underline{\alpha^i}}) \cap $supp$(\underline{x}^{\underline{\beta^i}}) = \emptyset$. The
Newton polytope of each toric hypersurface $h^i$ is a segment
$\tau^i$, connecting $\underline{\alpha^i}$ and
$\underline{\beta^i}$. Let $\tau^{*i}$ be the dual space in
$\mathbb{R}^{*d+p}_{\geq 0}$ to this segment and let $H^{*i}$ be the
hyperplane passing through $\tau^{*i}$. The equation of $H^{*i}$ is
$\sum_{k=1}^{d+p} (\alpha^i_k-\beta^i_k)x_k=0$. Set
$\tau^*:=\tau^{*1} \cap \cdots \cap \tau^{*p}$ and $H^*:=H^{*1} \cap
\cdots \cap H^{*p}$.
Let $N$ be the dual space to $M$ and let $\sigma$ be the dual cone
to $\check{\sigma}$. A primitive
element $n$ in $\sigma \cap N$ defines a discrete valuation $\nu$ of
$\mathbb{C}(V)$ by setting $\nu(\sum_{m \in F} a_m\chi^m)=\textnormal{min} \{ \langle m, n \rangle
\mid m \in F, a_m \neq 0\}$.

A finite set of valuations $\underline{\nu}$ in $\sigma$ induces a
Poincar\'{e} series $P^{\underline{\nu}}_V$ for $V$. On the other
hand valuations $\underline{\mu}$ in $\tau^*$ give rise to ambient
ideals $M(\underline{v}) \subset X=\mathbb{C}^{d+p}$ and hence to a Poincar\'{e} series $\mathcal{\mathcal{P}}^{\underline{\mu}}_{V}$
for $V$. We now show how both Poincar\'e series are related.
\begin{theorem} \label{thmlaatste}  % ${}$\\${}$     %\vspace*{-0.5cm}
%\begin{enumerate}
%\item
The cones $\sigma$ and $\tau^*$ are isomorphic
%\item
and under this isomorphism one has
$P^{\underline{\nu}}_V(\underline{t})=\mathcal{P}^{\underline{\mu}}_{V}(\underline{t})$.
%\end{enumerate}
\end{theorem}
\noindent \emph{Proof.} \quad
%\begin{enumerate}
%\item
We show that there is an isomorphism $\theta: \mathbb{R}^d
\rightarrow H^*$ that maps $\sigma$ to $\tau^*$. Let $\nu \in
\mathbb{R}^d$ and let $\mu=(\mu^1,\cdots,\mu^{d+p})$ be the vector
such that $\mu^k=\langle s_k, \nu \rangle$, $1 \leq k \leq d+p$.
Then obviously $\sum_{k=1}^{d+p} (\alpha^i_k-\beta^i_k)\mu^k=0$ for $1 \leq i \leq p$ and
this implies that $\mu \in H^*$. For the opposite direction, take
$\mu \in H^*$. The equality $\sum_{k=1}^{d+p}
(\alpha^i_k-\beta^i_k)\mu^k=0$ implies that there exists a $\nu \in
\mathbb{R}^d$ such that $\mu^k=\langle s_k,\nu \rangle$, $1 \leq k
\leq d+p$. This $\nu$ is then unique. As $\sigma$ is the dual
cone to $\check{\sigma}$, it follows that $\sigma$ maps to $\tau^*$.
%\item

The Poincar\'{e} series $\mathcal{P}^{\underline{\mu}}_{V}$ with
respect to the valuations $\mu_1,\cdots,\mu_r \in \tau^*$ is induced
by the ideals
\[J(\underline{v})=(\underline{x}^{\lambda} \mbox{ $|$ }
\langle \lambda, \mu_j\rangle \geq v_j, 1 \leq j \leq r) +
(h^1,\cdots,h^p).\] As $\langle \lambda, \mu_j \rangle = \langle
s,\nu_j \rangle$, where $s=\sum_{k=1}^{d+p} \lambda_k s_k$, it
follows that $\underline{x}^{\lambda} \in J(\underline{v})$ if and
only if $\varepsilon(\underline{x}^{\lambda})=\chi^s$, with $\langle
s, \nu_j \rangle \geq v_j$, $1 \leq j \leq r$.\\ The Poincar\'{e}
series $P^{\underline{\nu}}_V$ with respect to the corresponding valuations
$\nu_j$ is induced by the ideals
\[(\chi^s \mbox{ $|$ }
\langle s, \nu_j \rangle \geq v_j, 1 \leq j \leq r).\] It now
follows that both Poincar\'{e} series coincide. \hfill
$\blacksquare$
\\ ${}$
\newpage \begin{center}
\textsc{5. Embedded filtrations for nondegenerate singularities}
\end{center}
${}$ \\ %When in \cite{triomonodromy} it was discovered that the
%Poincar\'e series for an irreducible plane curve coincides with the
%zeta function of monodromy, the question arose if there exist other
%kinds of singularities for which a Poincar\'e series determines the
%topology of the singularity. This question has been answered
%positively for reducible plane curve singularities and for
%quasi-ordinary singularities. In \cite{pedro} one considers the
%Poincar\'e series with respect to the essential valuations over the
%singular locus and the essential valuations over the origin. Recall
%that the local ring for quasi-ordinary singularities (in particular
%irreducible plane curve singularities) is isomorphic to the local
%ring of a toric variety Spec $\mathbb{C}[S]$. In
%\cite{triomonodromy} and in \cite{pedro} one actually deduces the
%semigroup $S$ out of the Poincar\'e series. For these singularities,
%the semigroup $S$ yields the topology of the singularity.
%\\ \\
Let $h$ be the germ of a holomorphic function on $\C^d$ defining a hypersurface singularity $(V,o)$.
Let $\underline{\nu}=\{\nu_1,\ldots,\nu_r\}$ be the monomial valuations corresponding to the facets of the Newton polyhedron of $h$, including the non-compact facets. The centre of the valuations $\nu_i$ could then be a prime ideal different from the maximal ideal of $\mathcal{O}_{\C^d,o}$. We suppose that there is at least one compact facet such that there is at least one valuation with centre in the maximal ideal.
The definition of Poincar\'e series for affine toric varieties can be extended for a set of valuations which contains at least one valuation with centre in the maximal ideal. Indeed, notice that the $\chi(\mathbb{P}F_{\underline{v}})$ are then finite numbers (see also \cite{pedro} for an equivalent definition using graded rings) such that the Poincar\'e series $\sum_{\underline{v} \in
\mathbb{Z}^r}\chi(\mathbb{P}F_{\underline{v}})\underline{t}^{\underline{v}}$ is well-defined.
\begin{theorem}
Suppose that $h$ is nondegenerate with respect to
its Newton polyhedron $\mathcal{N}$ in the origin and that $\mathcal{N}$ has at least one compact facet. %We suppose that this Newton
%polyhedron contains at least one compact facet
Let $\underline{\nu}=\{\nu_1,\cdots,\nu_r\}$ be the monomial valuations on $\C^d$ induced by the facets of $\mathcal{N}$.
Then the Poincar\'e series $\mathcal{P}^{\underline{\nu}}_{V}(\underline{t})$ contains the same information as the Newton polyhedron of $h$ and in particular determines the zeta function
of monodromy of $h$.
\end{theorem}
\noindent \emph{Proof.} \quad Let $\nu_1$ be a valuation in $\underline{\nu}$ with centre in the maximal ideal.
The coefficient of $\underline{t}^{\underline{v}}$ in the series
\begin{eqnarray} \label{eqnpoincare}
\frac{(1-t_1^{q_1}\cdots t_r^{q_r})}{(1-t_1^{\nu_{1,1}}\cdots
t_r^{\nu_{r,1}})\cdots(1-t_1^{\nu_{1,d}}\cdots
t_r^{\nu_{r,d}})}\end{eqnarray} with $\underline{q}=\underline{\nu}(h)$,
can also be written as
\begin{eqnarray*}
\sum_{A \subset \{2,\cdots,r\}} (-1)^{\# A} \textnormal{dim
}\frac{M(\underline{v}-\underline{1}_A+
\underline{1})}{M(\underline{v}-\underline{1}_{A}+\underline{e_1}+ \underline{1})}-
\\ \sum_{A \subset \{2,\cdots,r\}} (-1)^{\# A} \textnormal{dim
}\frac{M(\underline{v}-\underline{q}-\underline{1}_A+
\underline{1})}{M(\underline{v}-\underline{q}-\underline{1}_{A}+\underline{e_1}+ \underline{1})}.
\end{eqnarray*}
Notice that these dimensions are finite because $\nu_1$ has centre in the maximal ideal.
One can now argue in the same way as in the proof of Theorem \ref{thmembedded} to obtain that
\[\mathcal{P}^{\underline{\nu}}_{V}(\underline{t})=
(1-\underline{t}^{\underline{q}})P^{\underline{\nu}}_{\C^d}(\underline{t}).\]
No factors cancel in (\ref{eqnpoincare}). Indeed, from the fact that $\underline{\nu}$ contains the valuations $(1,0,\ldots,0)$,\newline $\ldots,(0,\ldots,0,1)$ and at least one valuation where none of the entries is equal to $0$, one can deduce that some variable, let's say $x_1$, should divide $h$ and that the monomial $x_1$ should be contained in the support of $h$. This contradicts the hypothesis that the Newton polyhedron of $h$ contains at least one compact facet.
Picturing the hyperplanes
\[\nu_{j,1}x_1 + \cdots + \nu_{j,d}x_d = q_j, \qquad 1 \leq j \leq r,\]
the Newton polyhedron of $h$ is completely determined. In particular, for nondegenerate polynomials, one can compute the zeta
function of monodromy from the Newton polyhedron by the formula of Varchenko (see
\cite{varchenko}). \hfill $\blacksquare$

%\end{enumerate}

\footnotesize{

\end{document}